\documentclass[a4paper,reqno,12pt]{amsart}
\usepackage{amssymb}
\usepackage{brunnian}
\usepackage{enumitem}
\usepackage{float}
\restylefloat{figure,table}
\usepackage[colorlinks=true]{hyperref}
\hypersetup{allcolors=[rgb]{0.9,0.1,0.4}}
\usepackage{mathrsfs}
\usepackage[sort&compress,numbers]{natbib}
\usepackage[hang]{subfigure}
\usepackage{tikz}
\usetikzlibrary{arrows,%
  decorations.markings,%
  matrix,%
  snakes}
\usepackage{wrapfig}

\DeclareMathOperator{\bond}{Bond}
\DeclareMathOperator{\id}{id}
\DeclareMathOperator{\mor}{Mor}
\DeclareMathOperator{\nerve}{Nerve}
\DeclareMathOperator{\Pre}{Pre}
\DeclareMathOperator{\source}{source}
\DeclareMathOperator{\target}{target}
\let\leq=\leqslant
\let\geq=\geqslant

\newcommand{\C}{\mathscr{C}}
\newcommand{\diagram}[3]{\matrix (#1) [matrix of math nodes,row
  sep={#2},column sep={#3},text height=1.5ex,text depth=0.25ex]}
\newcommand{\F}{\mathscr{F}}
\renewcommand{\H}{\mathscr{H}}

\renewcommand{\P}{\mathscr{P}}
\renewcommand{\S}{\mathscr{S}}
\newcommand{\Sets}{\mathrm{Sets}}
\newcommand{\Sieve}{\mathcal{S}}
\newcommand{\sets}{\mathrm{Sets}}

\theoremstyle{definition}
\newtheorem*{definition}{Definition}
\newtheorem{example}{Example}

\title{On the Mathematics of Higher Structures}
\author{Nils A.\ Baas}
\address{Department of Mathematical Sciences, NTNU, N-7491 Trondheim,
  Norway}
\email{nils.baas@ntnu.no}
\date{April 23, 2019}

\begin{document}

\maketitle

\section{Introduction}
\label{sec:intro}

In a series of papers \cite{B1, B2, B3, B4, B5, B6, B7, B8, B9, B10,
  B11, B12, B13, B14, B15, Baas2018, BS2016} we have discussed higher
structures in science in general, and developed a framework called
Hyperstructures for describing and working with higher structures. In
\cite{B14} we discussed the philosophy behind higher structures and
formulated a principle in six stages --- the Hyperstructure Principle
--- for forming higher structures.

In this paper we will relate hyperstructures and the general principle
to known mathematical structures. We also discuss how they may give
rise to new mathematical structures and prepare a framework for a
mathematical theory.

Let us first recall from \cite{B14} what we think is the basic
principle in forming higher structures.

\section{The $\H$-Principle}
\label{sec:Hprin}

\begin{itemize}
\item[(I)] \emph{Observation and Detection.}\\

  \noindent Given a collection of objects that we want to study and
  give a structure. First we observe the objects and detect or assign
  their properties, states, etc.  This is the \emph{semantic} part of
  the process. Finally we may also select special objects.\\

\item[(II)] \emph{Binding.}\\

  \noindent A procedure to produce new objects from collections of old
  objects by ``binding'' them in some way. This is the
  \emph{syntactic} part of the process.\\
 
\item[(III)] \emph{Levels.}\\

  \noindent Iterating the described process in the following way:
  forming bonds of bonds and --- important! --- using the detected and
  observed properties at one level in forming the next level. This is
  iteration in a new context and not a recursive procedure. It
  combines syntax and semantics in forming a new level. Connections
  between levels are given by specifying how to dissolve a bond into
  lower level objects. When bonds have been formed to constitute a new
  level, observation and detection are like finding ``emergent
  properties'' of the process.\\
\end{itemize}

These three steps are the most important ones, but we include three
more in the general principle.

\begin{itemize}
\item[(IV)] \emph{Local to global.}\\

  \noindent Describing a procedure of how to move from the bottom
  (local) level through the intermediate levels to the top (global)
  level with respect to general properties and states. The importance
  of the level structure lies in the possibility of manipulating the
  systems levelwise in order to achieve a desired global goal or
  state. This can be done using ``globalizers'' --- an extension of
  sections in sheaves on \emph{Grothendieck sites} (see \cite{B2}).\\

\item[(V)] \emph{Composition.}\\

  \noindent A way to produce new bonds from old ones. This means that
  we can compose and produce new bonds on a given level, by ``gluing''
  (suitably interpreted) at lower levels. The rules may vary and be
  flexible due to the relevant context.\\

\item[(VI)] \emph{Installation.}\\

  \noindent Putting a level structure satisfying I--V on a set or
  collection of objects in order to perform an analysis, synthesis or
  construction in order to achieve a given goal. The objects to be
  studied may be introduced as bonds (top or bottom) in a level
  structure.\\

  \noindent\emph{Synthesis:} The given collection is embedded at the
  bottom level.\\

  \noindent\emph{Analysis:} The given collection is embedded at the
  top level.\\
\end{itemize}

Synthesis facilitates local to global processes and dually, analysis
facilitates global to local processes by defining localizers dual to
globalizers, see \cite{B1}.

The steps I--VI are the basic ingredients of what we call the
\emph{Hyperstructure Principle} or in short the
\emph{$\H$-principle}. (Corresponding to ``The General Principle'' in
\cite{B10}.)  In our opinion it reflects the basic way in which we
make or construct things.

Let us illustrate this in terms of category theory:

\begin{enumerate}
\item \emph{Observation} and \emph{detection}: we decide the structure
  of the objects like topological spaces, groups, etc.\\

\item \emph{Binding}: morphisms bind objects --- in an ordered way,
  continuous maps, homomorphisms, etc.\\

\item \emph{Levels}: we consider morphisms of morphisms of $\ldots$ in
  forming higher categories. Observation, detection and assignment
  become more indirect, but ought to play a more significant role.\\

\item \emph{Local to global}: at one level think of a Grothendieck
  sheaf on a site.\\

\item \emph{Composition}: composition of morphisms etc.\ in the
  ordinary sense.\\

\item \emph{Installation}: giving a collection of objects (like ``all
  groups'') a categorical structure.
\end{enumerate}

\section{A categorical implementation of the $\H$-principle}
\label{sec:categorical}

In order to illustrate how the $\H$-principle may be applied in an
ordinary categorical setting we take the following example from
\cite{B10}:

Let $\C$ be a category and $P: \C^{\text{op}} \to \sets$
a functor called a presheaf. The category of elements of $P$, denoted
by
\begin{equation*}
  \int_{\C} P \; \;,
\end{equation*}
is given as follows.
\begin{description}[style=nextline]
\item[Objects] $(C,p)$ where $C$ is an object in $\C$ and $p\in
  P(C)$.\\
\item[Morphisms] $(C',p') \to (C,p)$ are the morphisms
  $u \colon C' \to C$ in $\C$ such that $Pu \colon P(C) \to P(C')$ and
  $Pu(p) = p'$.
\end{description}
For this construction see \cite{MM}. 

Then a possible way to contruct a categorical hyperstructure is as 
follows: Start with a collection of objects $X_0$.
\begin{description}[style=nextline]
\item[Observation]
  \begin{align*}
    X_0 \rightsquigarrow  &\; \C_0\\
    &\; \text{$-$ category}\\
    &\\
    \C_0 \rightsquigarrow &\; \C_0 \text{ or}\\
    &\; \sets^{\C^{\text{op}}}\\
    &\; \C_0^J\\
    &\; \vdots\\
    &\\
    \Omega_0 \colon \C_o^{\text{op}} \to &\; \sets \text{ (Spaces,
      categories or other structures)}\\
    &\; \text{$-$ presheaf}
  \end{align*} 
\item[Binding]
  \begin{align*}
    \Gamma_0 = &\;\int_{\C_0} \Omega_0\\
    &\; \text{$-$ category of elements}\\
    &\\
    B_0 \colon \Gamma_0^{\text{op}} \to &\; \sets \text{
      (Spaces,$\ldots$)}\\
    &\; \text{$-$ presheaf}
  \end{align*}
\item[Levels]
  \begin{equation*}
    \C_1 = \int_{\Gamma_0} B_0
  \end{equation*}
  Iterating this process by making the appropriate choices we get a
  hyperstructure:
  \begin{equation*}
    \H = \{\C_0, \C_1,\ldots, \C_n \}
  \end{equation*}
  where
  \begin{align*}
    \C_m &= \int_{\Gamma_{m-1}} B_{m-1}\\
    &= \underset{\underset{\C_{m-1}}{\int} \Omega_{m-1}}{\int} B_{m - 1}
  \end{align*}
  for $1\leq m < n$.
\end{description}

In category theory it is often very useful to apply the nerve
construction to a category (even higher ones) in order to associate a
space from which topological information can be extracted. In the
present construction the ``nerve'' of $\H$ would mean the nerve of
$\C_n$ constructed inductively. The point to be made is that the
$\nerve(\H) = |\H|$ makes sense and may be useful in this context.

\section{From morphisms to bonds}
\label{sec:mor-bond}

In category theory we consider an ordered pair of objects $(X,Y)$ and
assign a set $\mor(X,Y)$ of morphisms. Intuitively the morphisms bind
the objects together. We suggest to extend the picture to a collection
of objects $\C = \{X_i\}_{i \in I}$.

The collection could be ordered or non-ordered. We prefer to present
here the ideas in the non-ordered case. Hence we assign a set of bonds
to the collection
\begin{equation*}
  B = B(\C).
\end{equation*}
$B$ may sometimes be empty.

The elements are mechanisms ``binding'' the collection in some way ---
extending morphisms. Let us look at some examples.

\subsection*{Relations}
A relation $R \subseteq X_1 \times \cdots \times X_n$ gives a bond of
tuples of elements $R(x_1,\ldots,x_n)$ if and only if
$(x_1,\ldots,x_n) \in R$.

\subsection*{Hypergraphs}
Here we are given a set of vertices and the edges are subsets of
vertices, and they serve as bonds of these vertices.

\subsection*{Subspaces}
Even more general, let $A_i$, $i = 1,\ldots,n$ be suitable subspaces of
$X$ and $A_i \subseteq X$. $X$ is then a bond of $\{A_i\}$. An
interesting case is when $A_i$ and $X$ are open subsets of a larger
space $Y$.

\subsection*{Simplicial complexes}
Given a simplicial complex $K$ based on vertices
$\{v_0,\ldots,v_q\}$. Then the simplices may be interpreted as
bonds.

\subsection*{Cobordisms}
Let $W, \{V_i\}_{i = 1,\ldots,k}$ be manifolds such that $\partial W =
\cup V_i$ ($V_i$ are the boundary components). We will then call $W$
a bond of $\{V_i\}$.

\subsection*{The basic idea:}
Instead of assigning a set $\mor(X,Y)$ to every ordered pair of
objects, we will assign a set of bonds to any collection of objects
--- finite, infinite or uncountable:
\begin{equation*}
  \bond(X,Y,Z,\ldots) \quad \text{or} \quad \bond(c\in \C)
\end{equation*}
$\C$ being a collection or parametrized family of objects. We may also
consider ordered collections or collections with other additional
properties. Bonds extend morphisms in categories and higher bonds
create levels and extend higher morphisms (natural transformations and
homotopies, etc.) in higher categories. This will be the basis for the
creation of new global states.

Bonds are more general than these examples. But prior to the bond
assignment is the process of observation, detection and assignment of
properties like: manifolds, subspaces, points, vertices, etc. This
will become more important when forming levels. Before studying level
formation we will discuss property and bond assignments.

Why do we need such an extension from graphs, higher categories,
etc. to hyperstructures? In previous papers \cite{B1, B2, B3, B4, B5,
  B6, B7, B8, B9, B10, B11, B12, B13, B14, B15, Baas2018, BS2016} ---
to which we refer the reader --- we have given many examples to
illustrate this: higher order links, higher cobordisms and many more
examples where we have group interactions instead of just pair
interactions. The essence is that many multiagent interactions require
a hyperstructure framework.

Here we just refer to these previous papers for examples and
motivation since our goal here is to discuss what we consider is the
essence of a philosophy of the mathematics of higher structures ---
outlining the possibilities for new constructions to be carried out in
the future.

\section{Property and bond assignments}
\label{sec:prop}

\subsection*{Properties}
By properties here we include: properties, states, phases,
etc. Collections we consider as subsets of some given set $X$, meaning
that a collection $S \in \P(X)$ --- the power set of $X$. In many
situations one may just consider structured subsets of $\P(X)$, but
the ideas remain the same. Similar to the example in Section
\ref{sec:categorical}. We may consider $\P(X)$ as a category with
inclusions as morphisms in some cases.

Even if the $\Omega$'s and $B$'s (to be defined later in this section)
are just general assignments we may ask how they behave with respect
to unions and intersections --- even if they are not functors. We may
look for analogues of pullback and pushout preservation. In many cases
we do not find this and it may lead to new kinds of mathematical
structures. This applies to both $\Omega$ and $B$ assignments.

We consider assignments
\begin{equation*}
  \Omega \colon \P(X) \to \sets
\end{equation*}
(or having target something more general like a higher
category). Should $\Omega$ be a functor, meaning that
\begin{equation*}
  S' \subseteq S
\end{equation*}
implies (contravariantly)
\begin{equation*}
  \Omega(S') \leftarrow \Omega(S)
\end{equation*}
or (covariantly)
\begin{equation*}
  \Omega(S') \to \Omega(S)?
\end{equation*}
In many situations this would be natural.\\

What about
\begin{equation*}
  \Omega(S' \cup S)
\end{equation*}
in terms of $\Omega(S')$ and $\Omega(S)$ where certainly $S' \cap S =
\emptyset$ is allowed?

\begin{figure}[H]
  \centering
  \begin{tikzpicture}[thick]
    \draw (0,0) circle(2cm and 1.5cm);
    \draw (2.5,0) circle(2cm and 2.5cm);
    \node at (-0.95,0.5){$S_1$};
    \node at (3,1.5){$S_2$};
    \draw[->] (-0.5,-2.5) node[below] {$S_1 \cap S_2$} to[out=90,in=250]
      (1,-0.5);
  \end{tikzpicture}
  \caption{Collections for state and bond assignments.}
  \label{fig:interection}
\end{figure}
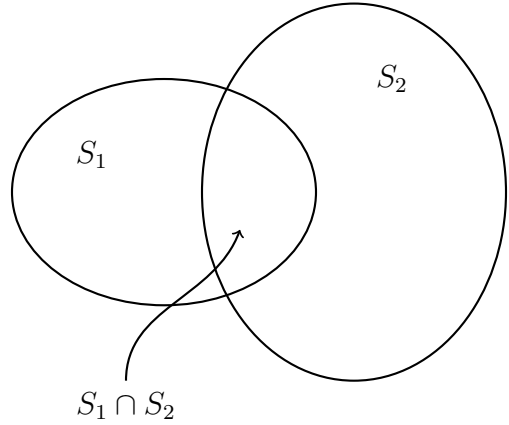

How does $\Omega(S_1 \cup S_2)$ relate to $\Omega(S_1)$, $\Omega(S_2)$
and $\Omega(S_1 \cap S_2)$?

\begin{enumerate}
\item\label{item1} If $\Omega$ is a covariant functor then
  \begin{center}
    \begin{tikzpicture}
      \begin{scope}[xshift=-3cm]
        \diagram{d}{1.5em}{1.5em}{
          S_1 \cap S_2 & S_1\\
          S_1 & S_1 \cup S_2\\
        };

        \path
          (d-1-1) edge[right hook->] (d-1-2)
                       edge[left hook->] (d-2-1)
          (d-1-2) edge[left hook->] (d-2-2)
          (d-2-1) edge[right hook->] (d-2-2);
      \end{scope}

      \node at (-0.25,0) {$\mapsto$};

      \begin{scope}[xshift=3cm]
        \diagram{dd}{1.5em}{1.5em}{
          \Omega(S_1 \cap S_2) & \Omega(S_1)\\
          \Omega(S_1) & \Omega(S_1 \cup S_2)\\
        };

        \path[->]
          (dd-1-1) edge (dd-1-2)
                         edge (dd-2-1)
          (dd-1-2) edge (dd-2-2)
          (dd-2-1) edge (dd-2-2);
      \end{scope}
    \end{tikzpicture}
  \end{center}
  (where $\Omega(S_1 \cup S_2) = \Omega(S_1) \sqcup_{\Omega(S_1 \cap
    S_2)} \Omega(S_2)$) which in some situations may be required to be
  a pushout.

  If $\Omega$ is contravariant, we may require a pullback:
  \begin{center}
    \begin{tikzpicture}
      \diagram{d}{1.5em}{1.5em}{
        \Omega(S_1 \cup S_2) & \Omega(S_2)\\
        \Omega(S_1) & \Omega(S_1 \cap S_2)\\
      };

      \path[->]
        (d-1-1) edge (d-1-2)
                     edge (d-2-1)
        (d-1-2) edge (d-2-2)
        (d-2-1) edge (d-2-2);
    \end{tikzpicture}
  \end{center}
  But there are situations in the general setting where none of these
  conditions are satisfied. We need to go beyond (co)-presheaves.
\item\label{item2} In some situations one may require a function or
  assignment $\varphi$ such that
  \begin{equation*}
    \Omega(S_1\cup S_2) = \varphi
    \bigl(\Omega(S_1),\Omega(S_2),\Omega(S_1\cap S_2)\bigl).
  \end{equation*}
  $\varphi$ may be thought of as a generalized limit in particular in
  the case of a union of an arbitrary collection of $S$'s.

  Properties or elements in $\Omega(S_1\cup S_2)$ not in or coming
  from $\Omega(S_1)$ or $\Omega(S_2)$ may be thought of as
  ``emergent'' properties.
\end{enumerate}

The theory should be developed in both the cases \ref{item1} and
\ref{item2}. In general the only assignment of ``emergent'' properties
is by ``observation'' --- ``the whole is more than the sum of its parts.''

\subsection*{Bonds}
We now consider collections $S$ with a property $\omega$, $\omega \in
\Omega(S)$ and form
\begin{equation*}
  \Gamma = \{ (S,\omega) \mid \omega \in \Omega(S)\}.
\end{equation*} 
We want to study the ``mechanisms'' that can bind the elements of $S$
together to some kind of unity. This is done by an assignment
\begin{equation*}
  B \colon \Gamma \to \sets,
\end{equation*}
where $B(S,\omega)$ is the set of bonds of $S$.

If $\Omega$ and $B$ are both functors we proceed by known mathematical
tools. If one of them or both fail to be functors we need to develop
new mathematical methods.

If
\begin{center}
  \begin{tikzpicture}
    \diagram{d}{1.5em}{1.5em}{
      (S_1 \cap S_2, \omega_{12}) & (S_1,\omega_1)\\
      (S_2,\omega_2) & (S_1 \cup S_2, \omega)\\
    };

    \path[->]
      (d-1-1) edge (d-1-2)
                   edge (d-2-1)
      (d-1-2) edge (d-2-2)
      (d-2-1) edge (d-2-2);
  \end{tikzpicture}
\end{center}
 it is sometimes natural to require that
\begin{center}
  \begin{tikzpicture}
    \diagram{d}{1.5em}{1.5em}{
      B(S_1 \cap S_2, \omega_{12}) & B(S_1,\omega_1)\\
      B(S_2,\omega_2) & B(S_1 \cup S_2, \omega)\\
    };

    \path[->]
      (d-1-1) edge (d-1-2)
                   edge (d-2-1)
      (d-1-2) edge (d-2-2)
      (d-2-1) edge (d-2-2);
  \end{tikzpicture}
\end{center}
is a pushout, or
\begin{center}
  \begin{tikzpicture}
    \diagram{d}{1.5em}{1.5em}{
      B(S_1 \cap S_2, \omega_{12}) & B(S_1,\omega_1)\\
      B(S_2,\omega_2) & B(S_1 \cup S_2, \omega)\\
    };

    \path[->]
      (d-1-2) edge (d-1-1)
      (d-2-1) edge (d-1-1)
      (d-2-2) edge (d-1-2)
                   edge (d-2-1);
  \end{tikzpicture}
\end{center}
a pullback. But sometimes these conventional notions fail and one may
proceed in different ways.

Bonds ($B$) (like morphisms) represent the syntactic part of the
structure. Observation ($\Omega$) --- missing in (Higher) Category
Theory --- represent the semantic part.

For property assignments $\Omega$ we may introduce operations: Given
$(S_1,\omega_1)$ and $(S_2,\omega_2)$, $\omega_1 \in \Omega(S_1)$ and
$\omega_2 \in \Omega(S_2)$, we may define
\begin{equation*}
  \omega_1 \circ \omega_2 = \varphi(\omega_1,\omega_2) \in
  \Omega(S_1\cup S_2), \quad \text{for $S_1$ and $S_2$ disjoint.}
\end{equation*}
Whenever a tensor product exists we may require:
\begin{equation*}
  \Omega(S_1\cup S_2) = \Omega(S_1) \otimes \Omega(S_2).
\end{equation*}
Whenever we introduce several levels properties will automatically
depend on previous properties in a cumulative way and take care of
levels. Bonds are different, composing and gluing at different
levels. Before elaborating that we need to discuss and specify the
formation of levels. First let us give two examples.

\begin{example}
  Given two sets of agents ($S_1$ and $S_2$) with specific skills (or
  products). In analogy with functorial assignments we will consider:
  \begin{enumerate}
  \item Let $\Omega$ assign collective skills. Then $\Omega(S_1)$ and
    $\Omega(S_2)$ will not necessarily map into $\Omega(S_1 \cap
    S_2)$. Hence no ``pullback property''.
  \item Let $\Omega$ assign individual skills to $S_1$ and $S_2$. Then
    $\Omega(S_1)$ and $\Omega(S_2)$ will not map into $\Omega(S_1 \cup
    S_2)$. Hence no ``pushout property''.
  \item Similarly for bonds, for example, formed by using skills to
    make certain products.
  \end{enumerate}
\end{example}

\begin{example}
  Given two sets of agents with specific skills (the $\Omega$-part)
  and a mechanism or organization binding them together to produce
  specific products (the $B$-part).

  The groups may intersect --- have agents in common --- but the
  intersection may be unable to produce the products. Hence no
  restriction maps or ``pullback property'' for bonds.

  Furthermore, we may consider the union of two groups which will
  clearly be able to produce the products of the groups, but the union
  may produce many more (for example composites). Hence, union is not
  preserved and no ``pushout property'' for bonds.
\end{example}

\section{Levels}
\label{sec:levels}

In higher categories we move from objects and morphisms to morphisms
of morphisms, etc. In the case of continuous maps we pass to
homotopies, homotopies of homotopies, etc. This is how higher levels
of structure arise.

In our situation we will now create higher levels by introducing bonds
of bonds, etc. Let us start with collections of objects from a basic
set $X_0$. Then we introduce as we described
\begin{equation*}
  \Omega_0, \Gamma_0, B_0.
\end{equation*}
We let the assignments --- whether functorial or not --- be sets, but
as we will point out later we may assign much more general
structures. (For example, $\infty$-groupoids or $\infty$-categories as
suggested by V.~Voevodsky in a private discussion.)

In forming the next level we define:
\begin{equation*}
  X_1 = \{b_0 \mid b_0 \in B_0(S_0,\omega_0), S_0\in \P(X_0) \text{
    and } \omega_0\in \Omega_0(S_0)\}.
\end{equation*}
Depending on the situation we now can choose $\Omega_1$ and $B_1$
according to what we want to construct or study and then repeat the
construction.

This is not a recursive procedure since new properties and bonds arise
at each level.

Hence a higher order architecture or structure of order $n$ is
described by:
\begin{equation*}
  \H_n \quad \colon \quad \begin{cases}
    X_0, \quad \Omega_0, \quad \Gamma_0, \quad B_0\\
    X_1, \quad \Omega_1, \quad \Gamma_1, \quad B_1\\[0.25cm]
    \hspace*{2cm} \vdots\\[0.25cm]
    X_n, \quad \Omega_n, \quad \Gamma_n, \quad B_n.
  \end{cases}
\end{equation*}
At the technical level we require that
\begin{equation*}
  B_i(S_i,\omega_i) \cap B_i (S_i',\omega_i') = \emptyset
\end{equation*}
for $S_i \neq S_i'$ (``a bond knows what it binds'') in order to
define the $\partial_i$'s below, or we could just require that the
$\partial_i$'s exist.

The level architectures are connected by ``boundary'' maps as follows:
\begin{equation*}
  \partial_i \colon X_{i + 1} \to \P(X_i)
\end{equation*}
defined by
\begin{equation*}
  \partial_i(b_i) = S_i \qquad \text{(dissolving bonds)}
\end{equation*}
and maps
\begin{equation*}
  I_i \colon X_i \to X_{i + 1}
\end{equation*}
such that $\partial_i \circ I_i = \id$.  $I_i$ gives a kind of
``identity bond''. $B_0$ may also contain identity bonds.

The extensions allowing bindings of subsets or subcollections of
higher power sets add many new types of architectures of
hyperstructures.  See \cite{B9,B2} for examples.

\begin{definition}
  We call the system
  \begin{equation*}
    \H_n = \{(X_i,\Omega_i,\Gamma_i,B_i,\partial_i) \mid i = 0,\ldots,n\}
  \end{equation*}
  a \emph{hyperstructure of order $n$}.
\end{definition}

This definition is made very general to illustrate the key idea. In
order to develop the definition and theory further mathematically
additional conditions will have to be added as pointed out in Section
\ref{sec:prop} and then it will branch off in several directions
depending on the situation under consideration, but with the
$\H$-structure as a common denominator. Our intention is also to cover
areas and problems outside of mathematics which again may give rise to
new mathematics.




\section{Composition of bonds}
\label{sec:comp}

In the study of collections of objects we emphasize the general notion
of bonds including relations, functions and morphisms.  We get richer
structures when we have composition rules of various types of bonds.
Such compositions should take into account the higher order
architecture giving bonds a level structure.

We experience this situation in higher categories where we want to
compose morphisms of any order.  Suppose that we are given two
$n$-morphisms $f$ and $g$.  They may not be compatible at level $n$
for composition in the sense that
\begin{equation*}
  \target(f) = \source(g).
\end{equation*}
But in a precise way we can iterate source and target maps to get down
to lower levels, and it may then happen that at level $p$ we have
\begin{equation*}
  \target_p^n(f) = \source_p^n(g).
\end{equation*}
Hence composition makes sense at level $p$ and we write the
composition rule as
\begin{equation*}
  \square_p^n
\end{equation*}
and the composed object as
\begin{equation*}
  f \, \square_p^n \, g.
\end{equation*}

In a similar way we can introduce composition rules for bonds in a
general hyperstructure $\H$.  Let $a_n$ and $b_n$ be bonds at level
$n$ in $\H$.  Then we get to the lower levels via the boundary maps
\begin{equation*}
  \partial_i \colon X_{i + 1} \to \P(X_i)
\end{equation*}
and search for compatibility in the sense that
\begin{equation*}
  \partial_p \circ \cdots \circ \partial_{n - 1} (a_n) = \partial_p
  \circ \cdots \circ \partial_{n - 1} (b_n)
\end{equation*}
or we may just require a weaker condition like
\begin{equation*}
  \partial_p \circ \cdots \circ \partial_{n - 1} (a_n)
  \cap \partial_p \circ \cdots \circ \partial_{n - 1} (b_n) \neq
  \emptyset
\end{equation*}
in order to have a composition defined:
\begin{equation*}
  a_n \, \square_p^n \, b_n
\end{equation*}

For bonds in a hyperstructure we may even compose bonds at different
levels: $a_m$, $b_n$ compatible at level $p$ via boundary maps, allow
us to define
\begin{equation*}
  a_m \, ^m\underset{p}{\square}^n \, b_n
\end{equation*}
as an $m$-bond for $m \geq n$.  Compositional rules are needed and
will appear elsewhere.

Composition may be thought of as a kind of geometric gluing.  We
consider the bonds as spaces, binding collections of families of
subspaces, these again being bonds, etc.  By the ``boundary'' maps we
go down to a level where these are compatible, gluable bond spaces
along which we may glue the bonds within the type of spaces we
consider. This applies for example to higher cobordisms.

Compositional rules will be needed, but they will depend on the
specific structures under study. For example we may require strict
associativity and/or commutativity or we may just require it up to a
higher bond. The point we are just trying to make is that there are a
lot of choices in the development of the further theory.

We have here for notational reasons suppressed the $\omega$'s
(properties/states), but they are included in a compatible way.

Therefore hyperstructures offer the framework for a new kind of higher
order gluing in which the level architecture plays a major role.  We
will pursue this in the next sections.

\section{States}
\label{sec:states}

Having introduced hyperstructures we may now assign states
(properties, etc.) to them:
\begin{equation*}
  \Lambda \colon \H \rightsquigarrow \S
\end{equation*}
where $\S$ is a structure representing the states --- in fact $\S$ may
be a level structure, a hyperstructure in itself. All assignments are
made level compatible. Furthermore, $\Lambda$ takes level to level and
may even be of a cumulative nature. The important point is assigning
states to bonds.

This means that
\begin{equation*}
  \Lambda = \{\Lambda_i\}, \qquad \S = \{\S_i\}
\end{equation*}
and
\begin{equation*}
  \begin{array}{r@{\ }c@{\ }l}
    \Lambda_0 & \text{takes values in} & \S_n\\
    & \vdots &\\
    \Lambda_i & \text{takes values in} & \S_{n - i}\\
    & \vdots &\\
    \Lambda_n & \text{takes values in} & \S_0.
  \end{array}
\end{equation*}
The degree of structure preservation may depend on the situation in
question.

Even if our starting hyperstructure $\H$ is very simple --- like a
multilevel decomposition of some space --- it may be very useful to
assign rather complex states in order to act on the system. This point
is dicussed in \cite[Section 5.1 --- $\H$-formation]{B14} where we
suggest that $\S$ may be a hyperstructure of higher types being
hyperstructures of hyperstructures $\ldots$

For state assignments there is a plethora of new possibilities,
extending assignments in topological quantum field theory (TQFT). In
such a level structure (hyperstructure) of states
\begin{equation*}
  \S = \{\S_0,\S_1,\ldots,\S_n\}
\end{equation*}
$\S_n$ represents the local states associated with the lowest level
bonds $B_0$, and $\S_0$ represents the global states associated with
the top bonds $B_n$.

As pointed out in \cite{B6,B2,B14} it is important to have level
connecting assignments making it possible to pass from local to global
states. Of course this is not always possible. We will discuss a way
of doing this by using generalized multilevel gluing. We use state
here in a general sense including observables and properties as
well. The important thing is that we in $\H$-structures have levels of
observables, states, properties, etc., not just local and global.

\section{Local to global}
\label{sec:locglob}

Hyperstructures are useful tools in passing from local situations to
global ones in collection of objects.  In this process the level
structure is important.  We will here elaborate the discussion of
multilevel state systems in \cite{B6} following \cite{B2}

In mathematics we often consider situations locally at open sets
covering a space and then glue together basically in one stroke ---
meaning there are just two levels local and global, no intermediate
levels.  In many situations dominated by a hyperstructure this is not
sufficient.  We need a more general hyperstructured way of passing
from local to global in general collections.

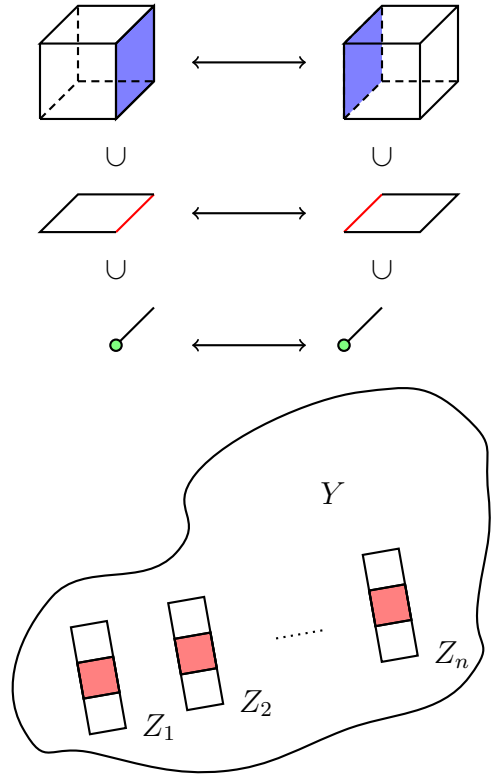
\begin{figure}[H]
  \centering
  \subfigure{
    \begin{tikzpicture}[thick]
      \begin{scope}
        \filldraw[fill=blue!50] (1,0) -- (1.5,0.5) -- (1.5,1.5) -- (1,1)
          -- cycle;
        \draw (0,0) -- (1,0) -- (1.5,0.5) -- (1.5,1.5) -- (0.5,1.5) --
          (0,1) -- (1,1) -- (1.5,1.5);
        \draw (0,0) -- (0,1);
        \draw (1,0) -- (1,1);
        \draw[densely dashed] (0,0) -- (0.5,0.5) -- (0.5,1.5);
        \draw[densely dashed] (0.5,0.5) -- (1.5,0.5);

        \draw[<->] (2,0.75) -- (3.5,0.75);

        \filldraw[fill=blue!50,densely dashed] (4,0) -- (4.5,0.5) --
          (4.5,1.5) -- (4,1) -- cycle;
        \draw (4,0) -- (5,0) -- (5.5,0.5) -- (5.5,1.5) -- (4.5,1.5) --
          (4,1) -- (5,1) -- (5.5,1.5);
        \draw (4,0) -- (4,1);
        \draw (5,1) -- (5,0);
        \draw[densely dashed] (4.5,0.5) -- (5.5,0.5);
      \end{scope}
  
      \begin{scope}[yshift=-0.5cm]
        \node[rotate=90] at (1,0){$\subset$};
        \node[rotate=90] at (4.5,0){$\subset$};
      \end{scope}
      
      \begin{scope}[yshift=-1.5cm]
        \draw[red] (1,0) -- (1.5,0.5);
        \draw (1,0) -- (0,0) -- (0.5,0.5) -- (1.5,0.5);

        \draw[<->] (2,0.25) -- (3.5,0.25);

        \draw[red] (4,0) -- (4.5,0.5);
        \draw (4,0) -- (5,0) -- (5.5,0.5) -- (4.5,0.5);
      \end{scope}

      \begin{scope}[yshift=-2cm]
        \node[rotate=90] at (1,0){$\subset$};
        \node[rotate=90] at (4.5,0){$\subset$};
      \end{scope}
      
      \begin{scope}[xshift=1cm,yshift=-3cm]
        \draw (0,0) -- (0.5,0.5);
        \filldraw[fill=green!50] (0,0) circle(0.075cm);

        \draw[<->] (1,0) -- (2.5,0);

        \draw (3,0) -- (3.5,0.5);
        \filldraw[fill=green!50] (3,0) circle(0.075cm);
      \end{scope}
    \end{tikzpicture}
  } \subfigure{
    \begin{tikzpicture}[scale=0.64,thick]
      \draw[rotate=15] (0.24396864,-1.4967076) .. controls
        (0.49101374,-2.4657116) and (0.665201,-1.8994577)
        .. (1.5439687,-2.3767076) .. controls (2.4227362,-2.8539574) and
        (3.017664,-3.2517743) .. (4.0039687,-3.4167075) .. controls
        (4.9902735,-3.581641) and (5.629609,-3.3427677)
        .. (6.6239686,-3.2367077) .. controls (7.6183286,-3.1306474) and
        (8.599862,-3.4573336) .. (9.383968,-2.8367076) .. controls
        (10.168076,-2.2160816) and (9.585438,-2.3502257)
        .. (9.943969,-1.4167076) .. controls (10.302499,-0.48318952) and
        (10.588963,-0.26462105) .. (10.743969,0.7232924) .. controls
        (10.898975,1.7112058) and (11.202379,1.6006097)
        .. (10.663969,2.4432924) .. controls (10.125558,3.2859752) and
        (9.734352,3.3049438) .. (8.743969,3.4432924) .. controls
        (7.753585,3.581641) and (5.849801,3.2754364)
        .. (5.0439687,2.6832924) .. controls (4.2381363,2.0911484) and
        (5.064207,1.7455099) .. (4.2239685,1.2032924) .. controls
        (3.3837304,0.661075) and (2.3561227,1.5378368)
        .. (1.5239687,0.9832924) .. controls (0.69181454,0.428748) and
        (-0.003076456,-0.5277036) .. (0.24396864,-1.4967076);

      \begin{scope}[xshift=2cm,yshift=-1.5cm,scale=0.75,rotate=10]
        \filldraw[fill=red!50] (0,1) -- (1,1) -- (1,2) -- (0,2) -- cycle;

        \draw (0,0) rectangle(1,3);
        \draw (0,1) rectangle(1,2);
        
        \node[right] at (1.125,0){$Z_1$};
      \end{scope}
      
      \begin{scope}[xshift=4cm,yshift=-1cm,scale=0.75,rotate=10]
        \filldraw[fill=red!50] (0,1) -- (1,1) -- (1,2) -- (0,2) -- cycle;

        \draw (0,0) rectangle(1,3);
        \draw (0,1) rectangle(1,2);
          
        \node[right] at (1.125,0){$Z_2$};
      \end{scope}
        
      \begin{scope}[xshift=6cm,yshift=-0.5cm,rotate=10]
        \draw[dotted] (0,1) -- (1,1);
      \end{scope}
      
      \begin{scope}[xshift=8cm,yshift=0cm,scale=0.75,rotate=10]
        \filldraw[fill=red!50] (0,1) -- (1,1) -- (1,2) -- (0,2) --  cycle;

        \draw (0,0) rectangle(1,3);
        \draw (0,1) rectangle(1,2);
        
        \node[right] at (1.125,0){$Z_n$};
      \end{scope}
        
      \node at (7,3.5){$Y$};
    \end{tikzpicture}
  }
  \caption{Gluing possibility at various levels}
  \label{fig:levels}
\end{figure}

Let us offer two of our intuitions regarding this process.
Geometrically we think of a multilevel nested family of spaces, like
manifolds with singularities represented by manifolds with multinested
boundaries or just like higher dimensional cubes with iterated
boundary structure (corners, edges,$\ldots$).  With two such
structures we may then glue at the various levels of the nesting
(Figure \ref{fig:levels}).
Furthermore, study how states and properties may be ``globalized'',
meaning putting local states coherently together to global states.

Biological systems are put together by multilevel structures from
cells into tissues, organs etc.\ constituting an organism.  Much of
biology is about understanding how cell-states determine organismic
states.  The hyperstructure concept is in fact inspired by biological
systems.

In order to extend the discussion of multilevel state systems in
\cite{B6} we need to generalize and formulate in a hyperstructure
context the following mathematical notions (see, for example,
\cite{MM}):
\begin{itemize}
\item Sieve
\item Grothendieck Topology
\item Site
\item Presheaf
\item Sheaf
\item Descent
\item Stack
\item Sheaf cohomology
\end{itemize}

Let us start with a given hyperstructure
\begin{align*}
  \H \quad \colon \quad & \{X_0,\ldots,X_n\}\\
  & \{\Omega_0,\ldots,\Omega_n\}\\
  & \{B_0,\ldots,B_n\}\\
  & \{\partial_0,\ldots,\partial_n\}\\
\end{align*}
We will now suggest a series of new definitions.

\begin{definition}
  \label{def:sieve}

  A \emph{sieve} on $\H$ is given as follows: at the lowest level
  $X_0$ a sieve $\Sieve$ on a bond $b_0 (=b_0(S_0,\omega_0))$ is given
  by families of bonds $\{ b_0^{j_0} \}$ (covering families) and
  $\beta_1$'s are compositional bonds in the family such that
  \begin{equation*}
    \beta_1(\{ b_0^{j_0} \}, b_0) \text{ --- the composition ---} 
  \end{equation*}
  is also in the family. $b_0$ may also be replaced by a family of
  bonds. ($b_0$ may also be an identity bond.)

  Bond composition with $\{b_0^{j_0}\}$ will produce new families in
  the sieve.

  A \emph{sieve on $\H$} is then a family of such sieves
  $(\Sieve_k)_{k = 1,\ldots,n}$ --- one for each level.
\end{definition}

We postpone connecting the levels until the definition of a
Grothendieck topology, but this could also have been added to the
sieve definition.

\begin{definition}
  \label{def:Grothendieck_topology}

  A \emph{Grothendieck topology} on $\H$ is given as follows: first we
  define a Grothendieck topology for each level of bonds.  Consider
  level $0$: to every bond $b_0$ we assign a collection of sieves
  $J(b_0)$ such that
  \begin{itemize}
    \item[(i)] (maximality), the maximal sieve on $b_0$ is in $J(b_0)$
    \item[(ii)] (stability), let $S \in J(b_0)$, $b_1(b_0',b_0)$, then in
      obvious notation
      \begin{equation*}
        b_1^\ast (S) \in J(b_0')
      \end{equation*}
    \item[(iii)] (transitivity), let $S \in J(b_0)$ and $R$ any sieve
      on $b_0,b_0'$ an element of a covering family in $S$, $b_1^\ast
      (R) \in J (b_0')$ for all $b_1$ with $b_1(b_0',b_0)$, then $R
      \in J(b_0)$.
  \end{itemize}
  We call $J(b_0)$ a $J$-covering of $b_0$.

  This gives a Grothendieck topology for all levels of bonds, and we
  connect them to a structure on all of $\H$ by defining in addition
  an assignment $J$ of $(b_0,\ldots,b_n)$ where $b_i \in \partial_i
  b_{i + 1}$.

  $J(b_0,\ldots,b_n)$ consists of families of sieves $\{b_0^{j_0}\}
  \in J(b_0),\ldots,\{b_n^{j_n}\}\in J(b_n)$ and bonds
  \begin{equation*}
    \beta_1,\ldots,\beta_{n + 1}
  \end{equation*}
  such that
  \begin{equation*}
    \beta_1(b_0,\{b_0^{j_0}\}),\ldots,\beta_{n +
      1}(b_n,\{b_n^{j_n}\})
  \end{equation*}
  and $b_i^{j_i} \in \partial_i b_{i + 1}^{j_{i + 1}}$.  In a diagram
  we have
  \begin{center}
    \begin{tikzpicture}[descr/.style={fill=white,inner sep=2.5pt}]
      \diagram{d}{2em}{2.5em}{
        b_n & b_{n - 1} & \cdots & b_0\\
        \{b_n^{j_n}\} & \{b_{n - 1}^{j_{n - 1}}\} & \cdots &
        \{b_0^{j_0}\}\\
        J(b_n) & J(b_{n - 1}) && J(b_0).\\
      };

      \path[->,midway,font=\scriptsize]
        (d-1-1) edge node[above]{$\partial$} (d-1-2)
                     edge[<->] node[right]{$\beta_{n + 1}$} (d-2-1)
        (d-1-2) edge node[above]{$\partial$} (d-1-3)
                     edge[<->] node[right]{$\beta_n$} (d-2-2)
        (d-1-3) edge node[above]{$\partial$} (d-1-4)
        (d-1-4) edge[<->] node[right]{$\beta_1$} (d-2-4)
        (d-2-1) edge node[above]{$\partial$} (d-2-2)
                     edge[-,white]
                     node[descr,sloped,text=black]{$\in$} (d-3-1)
        (d-2-2) edge node[above]{$\partial$} (d-2-3)
                     edge[-,white]
                     node[descr,sloped,text=black]{$\in$} (d-3-2)
        (d-2-3) edge node[above]{$\partial$} (d-2-4)
        (d-2-4) edge[-,white]
                     node[descr,sloped,text=black]{$\in$} (d-3-4);
    \end{tikzpicture}
  \end{center}
  Clearly there are many possible choices of Grothendieck topologies,
  and they will be useful in the gluing process and the creation of
  global states.  Examples will be discussed elsewhere, our main point
  here is to outline the general ideas.
\end{definition}

\begin{definition}
  \label{def:site}

  $(\H,J)$ is called a \emph{hyperstructure site} when
  $J$ is a Grothendieck topology on the hyperstructure $\H$.\sloppy
\end{definition}

Given
\begin{equation*}
  \S = \{ \S_0, \S_1,\ldots,\S_n\}
\end{equation*}
$\S_i$ being a hyperstructure and assignments such that
\begin{equation*}
  \begin{array}{r@{\ }c@{\ }l}
    \Lambda_0 & \text{takes values in} & \S_n\\
    & \vdots &\\
    \Lambda_i & \text{takes values in} & \S_{n - i}\\
    & \vdots &\\
    \Lambda_n & \text{takes values in} & \S_0.
  \end{array}
\end{equation*}
Sometimes we may also assume that $\S$ is organized into a
hyperstructure.  We assume that we have bond compatibility of the
$\Lambda_i$'s, preservation of bond composition and level connecting
assignments $\delta_i$ (``dual'' to the $\partial_i$'s and acting on
collections of bond ``states'') depending on the Grothendieck topology
$J$:
\begin{center}
  \begin{tikzpicture}
    \diagram{d}{2.5em}{2.5em}{
      \S_0 & \S_1 & \cdots & \S_n.\\
    };

    \path[->,midway,above,font=\scriptsize]
      (d-1-2) edge[decorate, decoration={snake, amplitude=.4mm,
        segment length=3mm, post length=1mm}] node{$\delta_1$}
        (d-1-1)
      (d-1-3) edge[decorate, decoration={snake, amplitude=.4mm,
        segment length=3mm, post length=1mm}] node{$\delta_2$}
        (d-1-2)
      (d-1-4) edge[decorate, decoration={snake, amplitude=.4mm,
        segment length=3mm, post length=1mm}] node{$\delta_n$}
        (d-1-3);
  \end{tikzpicture}
\end{center}
The $\delta_i$'s may be cumulative functional or relational
assignments, and the $\S_i$'s often have an algebraic structure.  In
the simplest case all the $\S_i$'s could just be $\Sets$. In defining
the $\delta$'s levels matter in a cumulative way and the $\delta$'s
may be seen as level connectors and regulators. See also \cite{B14}.

We consider the $\Lambda_i$'s as a kind of ``\emph{level presheaves}''
and the $\delta_i$'s giving a kind of ``\emph{global matching
  families}'' --- between levels in addition to levelwise matching.
However, if we have ``functional'' assignment connectors
$\hat{\delta}_i$'s on $\H$:
\begin{center}
  \begin{tikzpicture}
    \diagram{d}{2.5em}{2.5em}{
      \S_0 & \S_1 & \cdots & \S_n\\
    };

    \path[->,midway,above,font=\scriptsize]
      (d-1-2) edge node{$\hat{\delta}_1$} (d-1-1)
      (d-1-3) edge node{$\hat{\delta}_2$} (d-1-2)
      (d-1-4) edge node{$\hat{\delta}_n$} (d-1-3);
  \end{tikzpicture}
\end{center}
means that we get a unique state of global bond objects --- like an
\emph{amalgamation} for presheaves but here across levels in addition
to levelwise amalgamation.  Global bonds are ``covered'' as follows
(see \cite{B6})
\begin{center}
  \begin{tikzpicture}
    \diagram{d}{2.5em}{2.5em}{
      \{ b(i_n) \} & \{ b(i_{n - 1},i_n) \} & \cdots & \{
      b(i_0,\ldots,i_n) \}\\
    };

    \path[->,midway,above,font=\scriptsize]
      (d-1-1) edge node{$\partial_{n -1}$} (d-1-2)
      (d-1-2) edge node{$\partial_{n - 2}$} (d-1-3)
      (d-1-3) edge node{$\partial_0$} (d-1-4);
  \end{tikzpicture}
\end{center}
and states are being levelwise globalized in a cumulative way by
\begin{center}
  \begin{tikzpicture}
    \diagram{d}{2.5em}{2em}{
      \Lambda_n(\{ b(i_n) \}) & \Lambda_{n - 1}(\{ b(i_{n - 1}, i_n)
      \}) & \cdots & \Lambda_0 (\{ b(i_0,\ldots,i_n \}).\\
    };

    \path[->,midway,above,font=\scriptsize]
      (d-1-2) edge node{$\hat{\delta}_1$} (d-1-1)
      (d-1-3) edge node{$\hat{\delta}_2$} (d-1-2)
      (d-1-4) edge node{$\hat{\delta}_n$} (d-1-3);
  \end{tikzpicture}
\end{center}

With a slight abuse of notation we write this as
\begin{equation*}
  \Lambda \colon (\H,J) \to \S
\end{equation*}
and define $\Lambda = \{\Lambda_i\}$ as a \emph{``presheaf''} on
$(\H,J)$ ($\Pre(\H,J)$) and when
\begin{equation*}
  \Delta = \{\hat{\delta}_i\}
\end{equation*}
exists we have a unique global bond state.  This is like a
sheafification condition and we call $(\Delta,\Lambda)$ a
\emph{globalizer} of the site $(\H,J)$ with respect to $\Lambda$.

$\Lambda$ with $\Delta$ extends the sheaf notion here, gluing
within levels and between levels.

\emph{A globalizer is a kind of higher order or hyperstructured sheaf
  covering all the levels. Dually we may also introduce ``localizers''
in a similar way.}

The existence of $\Delta$ contains the global gluing data and hence
corresponds to what is often called \emph{descent conditions} and the
hyperstructure collection $\S$ extends the notion of a \emph{stack}
over $\H$. The ``internal'' $\Omega$-property assignments may also be
required to satisfy these globalizing conditions depending on the
situation, sometimes we omit them notationally. The details may be
worked out in several directions.

Topological quantum field theories are examples of this kind of
assignments. When higher cobordism categories of manifolds and
cobordisms with boundaries, e.g.\ cobordism categories with
singularities (see \cite{B6}), are considered the assignments may take
values in some ``algebraic'' higher category like higher vectorspaces
or higher factorization algebras.

Suppose that we have an assignment
\begin{equation*}
  \Lambda \colon \H \to \S
\end{equation*}
and consider a bond $b_i$ at level $i$ in $\H$:
\begin{equation*}
  \partial_i b_i = \{b_{i - 1}^j\} \quad \text{for all $i$}.
\end{equation*}
Then a globalizer will give an assignment
\begin{equation*}
  \prod_j \Lambda_{i - 1} (b_{i - 1}^j) \xrightarrow{\delta_{n - i +
      1}} \Lambda_i(b_i). 
\end{equation*}

This shows that from a family of ``things'' of one kind, one can make
a ``thing'' of another (higher) kind at a higher level. One may view
this as a vast generalization of the concept of an operad (see
\cite{Leinster}).

If the $\S_i$'s have a tensor type product we should require:
\begin{equation*}
  \bigotimes_j \Lambda_{i - 1}(b_{i - 1}^j) \to \Lambda_i (b_i).
\end{equation*}
Sometimes when it makes sense
\begin{equation*}
  \S_k = \S_{k - 1},
\end{equation*}
like often in field theory, we may have
\begin{equation*}
  \Lambda_i (b_i) \in \bigotimes_j \Lambda_{i - 1} (b_{i - 1}^j)
\end{equation*}
and
\begin{center}
  \begin{tikzpicture}
    \diagram{d}{2.5em}{5em}{
      b_i & \lambda_i\\
      \{b_{i - 1}^j\} & \{\lambda_{i - 1}^j\}\\
    };

    \path[font=\scriptsize]
      (d-1-1) edge[snake=coil,segment aspect=0,segment
        amplitude=0.5pt,->] node[midway,left]{$\partial$} (d-2-1)
      (d-1-1) edge[snake=coil,segment aspect=0,segment
        amplitude=0.5pt,->] node[midway,above]{$\Lambda_i$} (d-1-2)
      (d-2-2) edge[snake=coil,segment aspect=0,segment
        amplitude=0.5pt,->] node[midway,right]{$\delta$} (d-1-2)
      (d-2-1) edge[snake=coil,segment aspect=0,segment
        amplitude=0.5pt,->] node[midway,below]{$\Lambda_{i - 1}$}
        (d-2-2);
  \end{tikzpicture}
\end{center}
extending pairings in TQFTs.

Also the ``internal'' property and state assignments in a
hyperstructure may be considered as extended multilevel field theories
\begin{equation*}
  \Omega_k \colon B_k \to \S_{n - k}
\end{equation*}
where then $\omega_k \in \Omega_k(b_k)$ and collections
$\{(b_k,\omega_k)\}$ form the next level.

A generalized field theory in this sense
\begin{equation*}
  \Lambda \colon \H \to \S
\end{equation*}
may be conceived as a bond between the hyperstructures $\H$ and
$\S$. This picture may be extended to bonds of families of
$\H$-structures
\begin{equation*}
  B(\{\H_i\})
\end{equation*}
where the $\H_i$'s could be a suitable mixture of geometric,
topological and algebraic hyperstructures. 

\section{Remarks}
\label{sec:rem}

\subsection{Installation}
\label{sub:inst}

This means that we just have a set or collection of objects --- $X$
--- that we want to study and work with. This may be facilitated by
organizing $X$ into a hyperstructure $\H(X)$ as argued in previous
papers \cite{B13, B12, B10, B8, B9, B5, B6, B2, B1, B15, B14, B7, B11,
  B3, B4, Baas2018, BS2016}. This is analogues to the useful process
of organizing a collection of objects into a category. Then one may
put structure assignments on $\H(X)$ again
\begin{equation*}
  \Lambda \colon \H(X) \to \S
\end{equation*}
and iterate whenever needed.

\subsection{$\H$-algebras}
\label{sub:Halg}

In an $\H$-structure with bonds $\{B_0,B_1,\ldots,B_n\}$ we may define
operations or products of bonds by ``gluing.''  If $b_n$ and $b_n'$
are bonds in $B_n$ that are ``gluable'' at level $k$, then we ``glue''
them into a new bond $b_n \, \square_k^n \, b_n'$:
\begin{center}
  \begin{tikzpicture}
    \diagram{d}{2.5em}{5em}{
      b_n & b_n'\\
      b_k & b_k'\\
    };

    \path[font=\scriptsize]
      (d-1-1) edge[snake=coil,segment aspect=0,segment
      amplitude=0.5pt,->] node[midway,left]{$\partial \circ \cdots
        \circ \partial$} (d-2-1)
      (d-1-2) edge[snake=coil,segment aspect=0,segment
      amplitude=0.5pt,->] node[midway,right]{$\partial \circ \cdots
        \circ \partial$} (d-2-2)
      (d-2-1) edge[out=-30,in=-150,<->]
      node[midway,below,xshift=6.75em]{``gluable''
        (having similar parts to be identified).} (d-2-2);
  \end{tikzpicture}
\end{center}
$(\H,\{\square_k^n\})$ gives new forms of higher algebraic
structures. We have \emph{level operations} $\{\square_k^k\}$ and
\emph{interlevel operations} $\{\square_k^n\}$.

For geometric objects $X$ and $Y$ one may define a ``fusion'' product
\begin{equation*}
  X \, \square_\H \, Y
\end{equation*}
by using installed $\H$-structures on $\H(X), \H(Y)$ and
$\H(X \sqcup Y)$, see \cite{B2}.

As pointed out in the previous section if in an $\H$-structure we are
given a bond $b_k$ binding $\{b_{k - 1}^i\}$ the state assignments
will give levelwise assignments connected via a globalizer
\begin{equation*}
  \Lambda_{k - 1}(\{b_{k - 1}^i\}) \rightsquigarrow
  \Lambda_k(b_k).
\end{equation*}
The globalizers act as generalized pairings connecting levels. In some
cases like factorization algebras connecting local to global
observables they may be isomorphisms (in perturbative field theories),
see \cite{AFR,Ginot}, but not in general.

An \emph{$\H$-algebra} will be an $\H$-structure $\H$ with ``fusion''
operations $\square = \{\square_k^n\}$.  One may also add a
``globalizer'' (see \cite{B2}) and tensor-type products as just
described. The combination of a tensor product and a globalizer is a
kind of extension of a ``multilevel operad.''

\subsection{Hidden $\H$-structures}
\label{sub:hidden}

In addition to the examples mentioned in Section \ref{sec:mor-bond}
there are well-known interesting structures that may be viewed as
hyperstructures:\\
\begin{itemize}
\item[(1)] Higher categories in general with objects, morphisms,
  morphisms of morphisms ($2$-morphisms), etc., see, for example,
  Lurie \cite{L1,L2}. Globalizers and localizers extend to the ideas
  of (iterated) spans, cospans and local systems in higher categories,
  see, for example, Lurie \cite{L1} and Haugseng \cite{H}.\\
\item[(2)] Higher cobordisms, cobordisms with singularities ---
  cobordisms of cobordisms $\ldots$ with iterated structural
  boundaries, see \cite{B6,B2}. Observables may be states, tangential
  properties, cohomological properties,$\ldots$\\
\item[(3)] Syzygies and resolutions in homological algebra are
  examples of structures of higher relations, see
  \cite{MacLane}. Hilbert's syzygy theorem states that if $M$ is a
  finitely generated module over a polynomial ring in $n$ variables
  over a field, then it has a free resolution of length $\leq n$. In
  our language: there is an installment of a hyperstructure on $M$ of
  order $\leq n$.

  Geometrically we see this for example in Adams resolutions coming
  from a (co)-homology theory.\\
\item[(4)] Higher spaces may be built up gluing or linking together
  spaces using (co)-homologically detected properties. For example
  gluing two spaces through subspaces connected by a map or relation
  with certain (co)-homological properties. This process may be
  iterated using possibly new (co)-homology theories forming new
  levels and one gets spaces with hyperstructures. Hyperstructures
  offer a method of describing a plethora of new spaces needed in
  various situations. One may for example take families of general
  spaces, manifolds or simplicial complexes and organize them into
  suitable $\H$-structures giving $\H$-spaces, $\H$-manifolds and $\H$
  simplicial complexes combining syntax (combinatorics) and semantics
  ((co-)homology, homotopy, $\ldots$).
\end{itemize}

\subsection{$\H$-spaces}
\label{sub:Hspaces}

What is a space?  This is an old and interesting question.  We will
here add some higher (order) perspectives.  Often spaces are given by
open sets, metrics, etc.  They all give rise to bindings of points:
open sets, ``binding'' its points, distance binding points, etc.

In many contexts (of genes, neurons, links, subsets and subspaces,
$\ldots$) it seems more natural to specify the binding properties of
space by giving a hyperstructure --- even in addition to an already
existing ``space structure''.  In order to emphasize the binding
aspects of space we suggest that a useful notion of space should be
given by a set $X$ and a hyperstructure $\H$ on it.  Such a pair
$(X,\H)$ we will call an \emph{$\H$-space}.  It tells us how the
points or objects are bound together, see \cite{Baas2018} for an
example.

Clearly there may be many such hyperstructures on a set. They may all
be collected into a larger hyperstructure --- $\H^{\text{Total}}$ ---
which in a sense parametrizes the others. Ordinary topological spaces
will be of order $0$ with open sets as bonds. Through the bonds one
may now study the processes like fusion and fission in the space.

Our key idea is that ``spaces'' and ``hyperstructures'' are intimately
connected.

In neuroscience one studies ``space'' through various types of cells:
place-, grid-, border-, speed-cells,$\ldots$, see \cite{B15}.
All this spatial information should be put into the framework of a
$\H$-spaces with for example firing fields as basic bonds.  As
pointed out, the \emph{binding} problem fits naturally in here,
similarly ``cognitive'' and ``evolutionary'' spaces defined by
suitable hyperstructures.  Higher cognition should be described by
$\H$-spaces as well.

From a mathematical point of view simplicial complexes are also a kind
of hyperstructure based on the vertices and the simplices being
bonds. In a simplex all subsets of vertices are subsimplices. We have
discussed in \cite{B5,B7} that many bonds do not have this
property. For example a Brunnian bond is a bond of say $n$ elements in
such a way that $(n - 1)$ are not bound together. These can be
realized as Brunnian links of various orders, see \cite{B5,
  BS2016}. We may therefore suggest the following:

\begin{definition}
  A \emph{Brunnian complex} consists of
  \begin{itemize}
  \item[(i)] A set of vertices
  \item[(ii)] A family of subsets $\F$ --- the set of simplices, such
    that singletons are in $\F$ and so is $\emptyset$.
  \end{itemize}
\end{definition}

This means only certain subsets are simplices, not all of them as in
simplicial complexes.

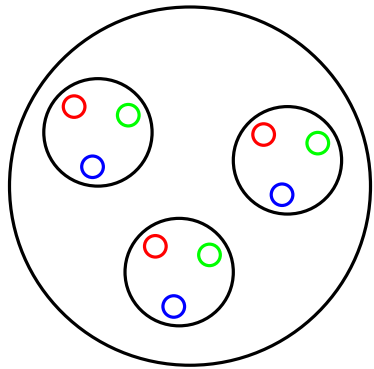
\begin{figure}[H]
  \centering
  \begin{tikzpicture}[scale=0.95]
    \draw[very thick] (0,0) circle(2.5cm);

    \begin{scope}[scale=0.3,xshift=-4.25cm,yshift=2.5cm]
      \draw[very thick] (0,0) circle(2.5cm);
      \draw[very thick, blue] (-0.25,-1.6) circle(0.5cm);
      \draw[very thick, red] (-1.1,1.2) circle(0.5cm);
      \draw[very thick, green] (1.4,0.8) circle(0.5cm);
    \end{scope}

    \begin{scope}[scale=0.3,xshift=4.5cm,yshift=1.2cm]
      \draw[very thick] (0,0) circle(2.5cm);
      \draw[very thick, blue] (-0.25,-1.6) circle(0.5cm);
      \draw[very thick, red] (-1.1,1.2) circle(0.5cm);
      \draw[very thick, green] (1.4,0.8) circle(0.5cm);
    \end{scope}
    
    \begin{scope}[scale=0.3,xshift=-0.5cm,yshift=-4cm]
      \draw[very thick] (0,0) circle(2.5cm);
      \draw[very thick, blue] (-0.25,-1.6) circle(0.5cm);
      \draw[very thick, red] (-1.1,1.2) circle(0.5cm);
      \draw[very thick, green] (1.4,0.8) circle(0.5cm);
    \end{scope}
  \end{tikzpicture}
  \caption{A Brunnian complex.}
  \label{fig:brunnian_complex}
\end{figure}

\begin{figure}[H]
  \centering
  \subfigure[Brunnian rings]{
    \begin{tikzpicture}[every path/.style={knot},scale=0.75]
      \pgfmathsetmacro{\brscale}{1.8}
      \foreach \brk in {1,2,3} {
        \begin{scope}[rotate=\brk * 120 - 120]
          \colorlet{chain}{ring\brk}
          \brunnianlink{\brscale}{120}
        \end{scope}
      }
    \end{tikzpicture}
    \label{fig:1stBrunnian}
  } \subfigure[2nd order Brunnian rings]{
    \includegraphics[width=0.35\linewidth]{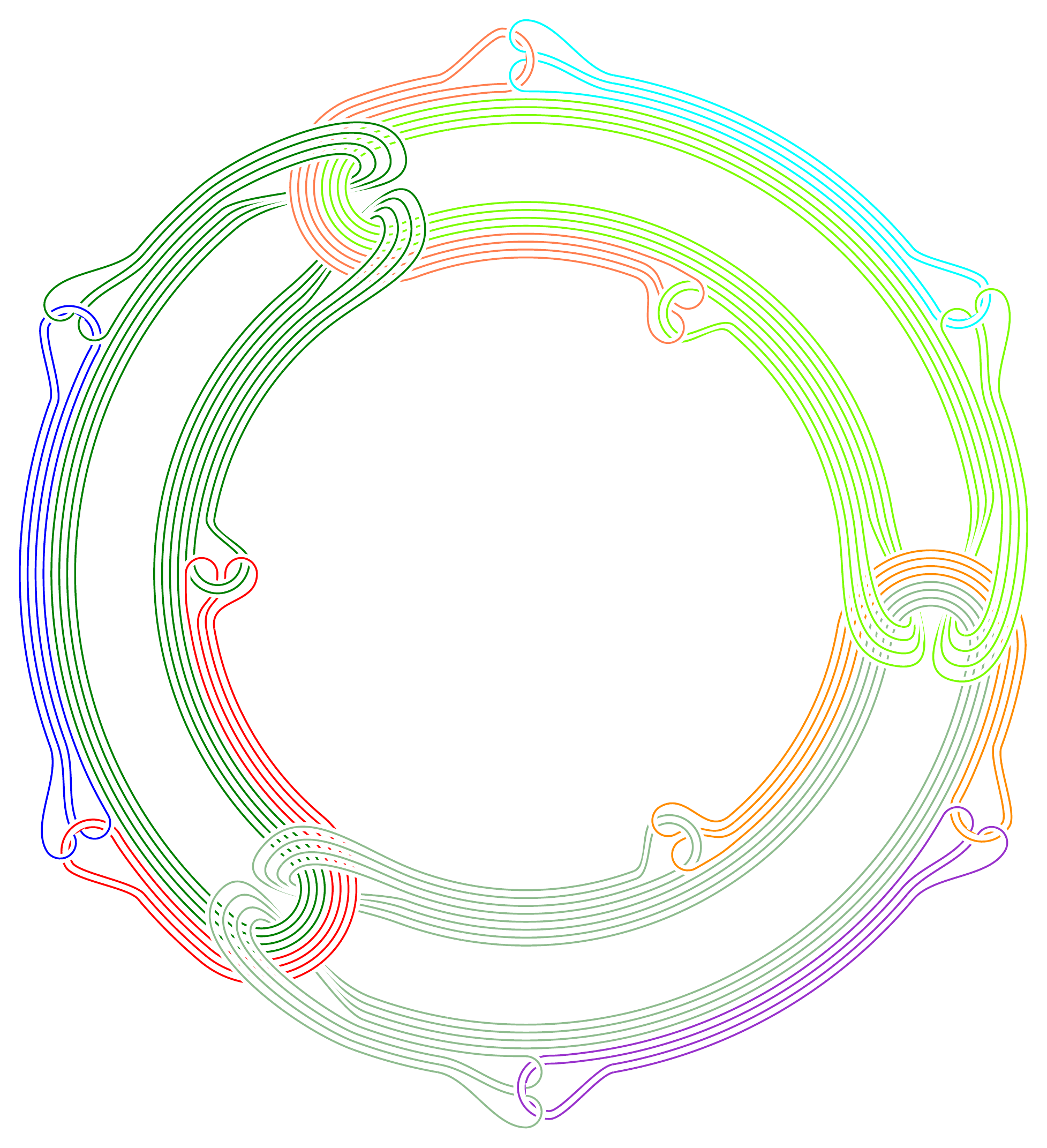}
    \label{fig:2ndBrunnian}
  }
  \caption{Links}
  \label{fig:links}
\end{figure}

In Figure \ref{fig:brunnian_complex} we have a $2$nd order Brunnian
complex of $9$ vertices and $3$ simplices, see Figure \ref{fig:links}
for the corresponding links.

\section{Conclusion}
\label{sec:conc}

The purpose of this paper is to introduce and formulate the basic
principles of higher structures occuring in science and nature in
general and in mathematics in particular. This suggests extensions of
known mathematical theory, but also leads to situations where new
mathematical theory has to be developed. This program of
Hyperstructures may go in many directions and we just consider this
paper as an eye opener of where to go in the future.

\subsection*{Acknowledgements}
I would like to thank P.~Cohen, D.~Sullivan and V.~Voevodsky for
interesting discussions at various stages of the development of the
mathematical aspects of the hyperstructure concept.

I would also like to thank M.~Thaule for his kind technical assistance in
preparing the manuscript. I would like to thank A.~Stacey for
producing Figure \ref{fig:links}.

\section*{Notes on the contributor}
\begin{wrapfigure}[8]{l}{3cm}
   \centering
   \vspace*{-0.5cm}
   \includegraphics[width=3cm]{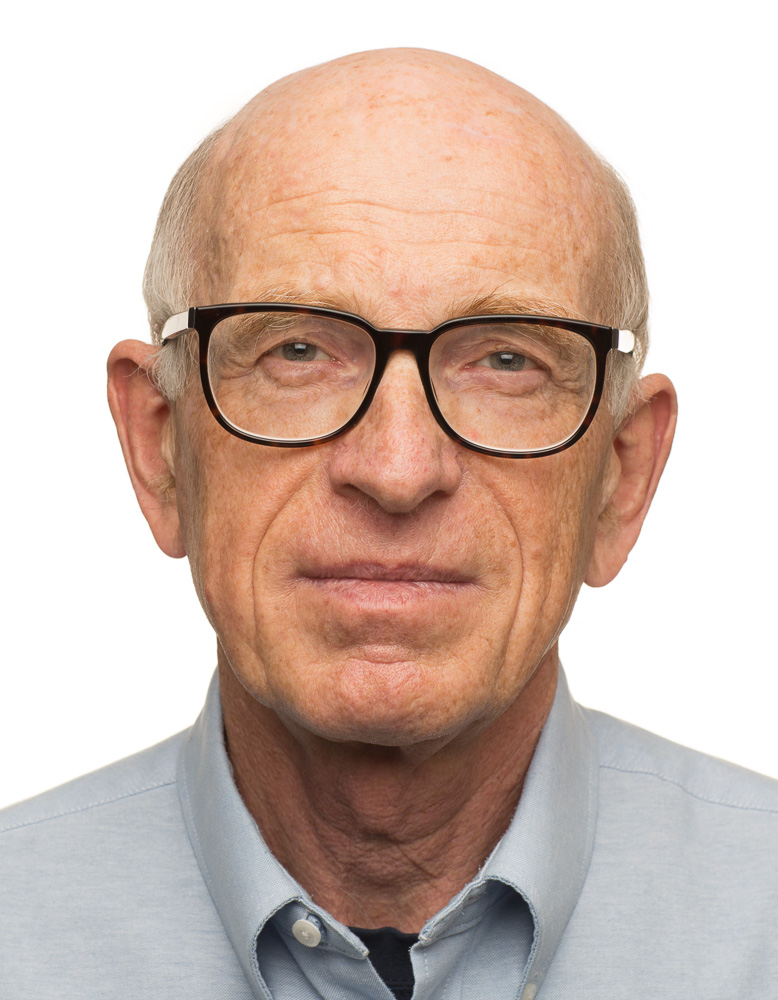}
\end{wrapfigure}

\noindent \textbf{Nils A.\ Baas} was born in Arendal, Norway, 1946.
He was educated at the University of Oslo where he got his final
degree in 1969.  Later on he studied in Aarhus and Manchester.  He was
a Visiting Assistant Professor at U. Va. Charlottesville, USA in
1971--1972.  Member of IAS, Princeton in 1972--1975 and IHES, Paris in
1975.  Associate Professor at the University of Trondheim, Norway in
1975--1977 and since 1977, Professor at the same university till date.
He conducted research visits to Berkeley in 1982--1983 and 1989--1990;
Los Alamos in 1996; Cambridge, UK in 1997, Aarhus in 2001 and 2004.
He was Member IAS, Princeton 2007, 2010, 2013 and 2016.  His research
interests include: algebraic topology, higher categories and
hyperstructures and topological data analysis.


%
%

\end{document}